\documentclass[12pt,oneside]{article}
\usepackage{amsmath,amssymb,amsfonts,amsthm}
\def\goth{\mathfrak}

\newcommand{\T}{{\cal T}}

\newcommand{\Real}{\mathbb R}

\newcommand{\To}{\longrightarrow}

\newcommand{\prof}{\noindent \textit{\textbf{Proof.\:}}}
\newcommand{\tm}{\T M}
\newcommand{\p}{\pi^{-1}(TM)}
\newcommand {\cp}{\mathfrak{X}(\pi (M))}
\def\o#1{\overline{#1}}

\def\ti#1{\tilde{#1}}

\def\x{{\goth X}(T_2M)}

\def\Section#1{\vspace{30truept}\addtocounter{section}{1}\setcounter{thm}{0}{\noindent\Large\bf
    \arabic{section}.~~#1}\par \vspace{12pt}}

\newtheorem{thm}{Theorem}[section]
\newtheorem{cor}[thm]{Corollary}
\newtheorem{lem}[thm]{Lemma}
\newtheorem{prop}[thm]{Proposition}

\title{ON GENERALIZED RANDERS MANIFOLDS}
\author{Aly A. Tamim\and Nabil L. Youssef}
\date{}

\begin{document}               
\bibliographystyle{plain}
\maketitle                     
\vspace{-1.45cm}
\begin{center}
{Department of Mathematics, Faculty of Science,\\ Cairo University,
Giza, Egypt.}
\end{center}

\begin{center}
Dedicated to Professor Radu Miron\\
on the occasion of his 70th birthday
\end{center}

\vspace{0.5cm} \maketitle
\smallskip


\vspace{13truept}\centerline{\Large\bf{Introduction}}\vspace{12pt}
\par
By a Randers' structure on a manifold $M$ we mean a Finsler
structure $L^*=L+\alpha$, where $L$ is a Riemannian structure and
$\alpha$ is a $1$-form on $M$. This structure was first introduced
by Randers ~\cite{[8]} from the standpoint of general relativity and
was investigated by several authors (~\cite{[2]}, ~\cite{[3]},
~\cite{[9]}, ...etc.) from the geometrical viewpoint. Numata
~\cite{[7]} studied Randers manifolds in the case where $L$ is a
locally Minkowskian structure on $M$.
\par
In this paper, we replace $L$ by a Finsler structure, calling the
resulting manifold a generalized Randers manifold. Such a manifold
was studied (using local coordinates) by Matsumoto ~\cite{[3]},
Tamim ~\cite{[10]} and Miron ~\cite{[6]}. Our aim is twofold. On one
hand, to pursue and develop in depth one of the present authors'
study ~\cite{[10]} of generalized Randers manifolds. On the other
hand, to apply the results obtained in a foregoing paper
~\cite{[12]} to generalized Randers manifolds to obtain some new
results in that domain. Among many results, we establish a necessary
and sufficient condition for a generalized Randers manifold to be a
general Landsberg manifold.
\par
It should be noticed that our approach is a global one. That is, it
does not make use of the local coordinate techniques (apart from the
proof of Theorem \ref{thm.1}).

\Section{Notations and Preliminaries}
\par
In this section, we give a brief account of the basic concepts
necessary for this work. For more details, refer to ~\cite{[1]} or
~\cite{[11]}. The following notations will be used throughout the
paper:

\begin{description}
    \item[] $M$: a differentiable manifold of finite dimension and of
class $C^\infty$.

    \item[]$\pi_M:TM \To M$: the tangent bundle of $M$.

    \item[] $\pi:\T M \To M$: the subbundle of nonzero
     vectors tangent to $M$.

    \item[]$P:\pi^{-1}(TM)\To \T M$: the
bundle, with base
     space $\T M$, induced by $\pi$ and $TM$

    \item[] $ \mathfrak{F}(M)$: the $\Real$-algebra of differentiable functions
     on $M$.

    \item[]$\mathfrak{X}(M)$: the $\mathfrak{F}(M)$-module of vector fields on
$M$.

    \item[] $\mathfrak{X}(\pi (M))$: the $\mathfrak{F}(\T M)$-module of differentiable
sections of $\pi^{-1}(TM)$.

\end{description}

Elements of $\mathfrak{X}(\pi (M))$ will be called $\pi$-vector
fields and will be denoted by barred letters $\overline{X}$. Tensor
fields on $\p$ will be called $\pi$-tensor fields. The fundamental
vector field is the $\pi$-vector field $\o\eta$ defined by
$\overline{\eta}(u)=(u,u)$ for all $u \in \T M $. The lift to $\p$
of a vector field $X$ on $M$ is the $\pi$-vector field $\o X$
defined by $ \overline{X}(u)=(u,X(\pi(u)))$.
\par
The vector bundles $\T M$ and $\p$ are related by the short exact
sequence
$$0\longrightarrow\p\stackrel{\gamma}\longrightarrow T(\tm)\stackrel{\rho}
\longrightarrow\p\longrightarrow0,$$ where the vector bundle
morphisms are defined by $\rho=(\pi_{\tm},d\pi)$ and
$\gamma(u,v)=j_u(v)$, where $j_u$ is the natural isomorphism
$j_u:T_{\pi_M(v)}M\longrightarrow T_u(T_{\pi_M(v)}M)$.
\par
Let $\nabla$ be an affine connection (or simply a connection) in the
vector bundle $\p$. We associate to $\nabla$ the map
$$K:\tm\longrightarrow\p:X\longmapsto\nabla_X\o\eta,$$
called the connection map of $\nabla$. A tangent vector $X\in
T_u(\tm)$ is said to be horizontal if $\ K(X)=0$. The connection
$\nabla$ is said to be regular if
$$T_u(\tm)=V_u(\tm)\oplus H_u(\tm)\qquad\forall u\in\tm,$$
where $V_u(\tm)$ and $H_u(\tm)$ are respectively the vertical and
horizontal spaces at $u$. If $M$ is endowed with a regular
connection, we can define a section $\beta$ of the morphism $\rho$
by $\ \beta=(\rho\mid_{H(\tm)})^{-1}$. It is clear that $\
\rho\circ\beta\ $ is the identity map on $\p$ and $\ \beta\circ\rho\
$ is the identity map on $H(\tm)$. Let ${\bf T}$ be the torsion form
and ${\bf R}$ the curvature transformation of the connection
$\nabla$. The horizontal and mixed torsion tensors, denoted
respectively by $S$ and $T$, are defined, for all $\o X,\o Y\in\cp$,
by:
$$S(\o X,\o Y)={\bf T}(\beta\o X,\beta\o Y),\quad
T(\o X,\o Y)={\bf T}(\gamma\o X,\beta\o Y).$$ The horizontal, mixed
and vertical curvature tensors, denoted respectively by $R$, $P$ and
$Q$, are defined, for all $\o X, \o Y,\o Z\in\cp$, by:
$$R(\o X,\o Y)\o Z={\bf R}(\beta\o X,\beta\o Y)\o Z,\ \
P(\o X,\o Y)\o Z={\bf R}(\gamma\o X,\beta\o Y)\o Z,\ \ Q(\o X,\o
Y)\o Z={\bf R}(\gamma\o X,\gamma\o Y)\o Z.$$
\par
If $c:I\longrightarrow M$ is a regular curve in $M$, its canonical
lift to $\tm$ is the curve $\ti c$ defined by $\ti c:t\longmapsto
dc/dt$. If $c$ is a geodesic in $M$, we shall denote by $\o V$ the
restriction of $\o\eta$ on $\ti c(t)$: $\o V=\o\eta\mid_{\ti c(t)}$.

\Section{Generalized Randers Manifolds}
\par
Let $(M,L)$ be a Finsler manifold. Let $g$ be the Finsler metric
associated with $(M,L)$. Using the bundle morphism $\gamma$, we
define the $\pi$-form:
\begin{equation}\label{eq.1}
\ell=dL\circ\gamma.
\end{equation}
One can easily show, for all $\o X\in\cp$, that
\begin{equation}\label{eq.2}
\ell(\o X)=L^{-1}g(\o X,\o\eta).
\end{equation}
The angular metric tensor $h$ is defined by:
\begin{equation}\label{eq.3}
h=g-\ell\otimes\ell.
\end{equation}
\par
Let $\nabla$ denote the Cartan's connection with respect to $g$.

\begin{lem}\label{lem.1}
For every $\ \o X,\o Y,\o Z\in\cp$, we have\:
\begin{description}
    \item[(a)]$\nabla_{\gamma\o X}L=\ell(\o X)$.

    \item[(b)] $(\nabla_{\gamma\o
X}\ell)(\o Y)=L^{-1}h(\o X,\o Y)$.

   \item[(c)]$(\nabla_{\gamma\o
X}h)(\o Y,\o Z)=-L^{-1}\{h(\o X,\o Y)\ell(\o Z)
       +h(\o X,\o Z)\ell(\o Y)\}$.

    \item[(d)] $\nabla_{\beta\o X}L=\nabla_{\beta\o X}\ell=\nabla_{\beta\o
X}h=0$.
\end{description}
\end{lem}

\par
Let $\delta$ be a given $1$-form on $M$. Let $\o b$ be the
$\pi$-vector field defined in terms of $\delta$ by:
$$\delta(X)=g(\o b,\o X)\qquad\forall X\in\x(M).$$
Writing
\begin{equation}\label{eq.4}
L^*=L+\alpha,\mbox{  where}\quad\alpha=g(\o b,\o\eta),
\end{equation}
$L^*$ defines a new Finsler structure (~\cite{[10]} and ~\cite{[6]})
on the manifold $M$. The Finsler manifold $(M,L^*)$ is called a
generalized Randers manifold and $(M,L)$ its associated Finsler
manifold.
\par
Using equations (\ref{eq.1})--(\ref{eq.4}), the $\pi$-tensors
$\ell$, $h$ and $g$ associated with $(M,L)$ and the corresponding
$\pi$-tensors associated with $(M,L^*)$ are related by:
\begin{equation}\label{eq.5}
  \left.
    \begin{array}{rcl}
         \ell^* &=& \ell+\omega,\mbox{  where  } \omega=d\alpha\circ\gamma\\
         h^*    &=& \tau h,\mbox{  where  } \tau=L^*L^{-1}\\
         g^*    &=& \tau(g-\ell\otimes\ell)+\ell^*\otimes \ell^*
    \end{array}
  \right\}
\end{equation}
\par

\begin{prop}\label{prop.1}
Let $\o m$ be the $\pi$-vector field defined by $\ \o m=\o
b-(\alpha/L^2)\o\eta,\ $ let $\nu$ be the $\pi$-form associated with
$\o m$ under the duality defined by the metric $\ g\ $ and let
$\phi$ be the $\pi$-form defined by $\
\phi=I-L^{-1}\ell\otimes\o\eta.$ Then, we have\:
\begin{description}
\item[(a)] $\ell(\o m)=0.$
\item[(b)] $\ell^*(\o m)=b^2-(\alpha/L)^2,\mbox{ where }b^2=g(\o b,\o b)$.
\item[(c)]$\nu(\o m)=b^2-(\alpha/L)^2.$
\item[(d)] $\phi(\o m)=\o m$.
\item[(e)]$\nu(\o X)=L\nabla_{\gamma\o X}\tau$.
\item[(f)]$\phi^*=\phi-{L^*}^{-1}\nu\otimes\o\eta$.
 \end{description}
\end{prop}

\par
For every $X\in\mathfrak{X}(TM)$ and $\o Y\in\cp$, let us write
\begin{equation}\label{eq.6}
{\nabla^*}_X\o Y=\nabla_X\o Y+U(X,\o Y),
\end{equation}
where $U$ is an ${\goth F}(\tm)$-bilinear mapping
$\mathfrak{X}(\tm)\times\cp \longrightarrow\cp$ representing the
difference between the two connections $\nabla^*$ and
$\nabla$.\newline For every $\o X,\o Y\in\cp$, we set
\begin{equation}\label{eq.7}
  \left.
      \begin{array}{rclcl}
A(\o X,\o Y) &:=& U(\gamma\o X,\o Y),\qquad B(\o X,\o Y) &:=& U(\beta\o X,\o Y)\\
     N(\o X) &:=& B(\o X,\o\eta),\qquad\qquad\quad\  N_0 &:=& N(\o\eta)
      \end{array}
  \right\}
\end{equation}
As a vector field $X$ on $\tm$ can be represented by $\ X=\gamma K
X+\beta\rho X$, it follows from (7) that
\begin{equation}\label{eq.8}
U(X,\o Y)=A(K(X),\o Y)+B(\rho X,\o Y),\qquad\forall
X\in{\mathfrak{X}(TM)},\ \o Y\in\cp.
\end{equation}
\par
From equations (\ref{eq.6}) and (\ref{eq.8}), we have

\begin{lem}\label{lem.2}
The two connections $\nabla$ and $\nabla^*$ are related by\:
$${\nabla^*}_X\o Y=\nabla_X\o Y+A(K X,\o Y)+B(\rho X,\o Y),$$
for all $X\in\mathfrak{X}(TM)$, $\o Y\in\cp$ .
\end{lem}
The $\pi$-tensor fields $A$ and $B$ will be determined explicitely
later.
\par
One can easily show that
\begin{equation}\label{eq.9}
\beta^*=\beta-\gamma\circ N.
\end{equation}
\par
Taking the definition of the torsion tensor $T$ into account and
using equations (\ref{eq.6}), (\ref{eq.9}) and (\ref{eq.5}), we get

\begin{prop}\label{prop.2}
For every $\o X,\o Y,\o Z\in\cp$, we have
\begin{description}

\item[(a)]$T^*(\o X,\o
Y)=T(\o X,\o Y)+A(\o X,\o Y)$.

\item[(b)] $T^*(N(\o X),\o Y)-T^*(N(\o
Y),\o X)=B(\o X,\o Y)-B(\o Y,\o X)$.

\item[(c)]$T^*(\o X,\o Y,\o Z)=\tau T(\o X,\o Y,\o Z)+\omega(T(\o X,\o Y))
\ell^*(\o Z)+A^*(\o X,\o Y,\o Z),$\newline where $T(\o X,\o Y,\o
Z):=g(T(\o X,\o Y),\o Z)$, $\ T^*(\o X,\o Y,\o Z):=g^*(T^*(\o X,\o
Y),\o Z)$ and \newline $A^*(\o X,\o Y,\o Z):=g^*(A(\o X,\o Y),\o Z)\
$.
\end{description}
\end{prop}

\begin{cor}\label{cor.1}
For all $\ \o X,\o Y,\o Z\in\cp$, we have:
\begin{description}

\item[(a)]$A^*(\o X,\o Y,\o Z)=A^*(\o X,\o Z,\o Y) +\omega(T(\o X,\o
Z))\ell^*(\o Y)-\omega(T(\o X,\o Y))\ell^*(\o Z).$

\item[(a)] $A^*(\o
X,\o Y,\o\eta)=-L^*\omega(T(\o X,\o Y)).$
\end{description}
\end{cor}
\par
Concerning the curvature tensors, we have

\begin{prop}\label{prop.3}
For every $\ \o X,\o Y,\o Z\in\cp$, we have
\begin{description}
\item[(a)] $R^*(\o X,\o Y)\o Z+P^*(N(\o X),\o Y)\o Z-P^*(N(\o
Y),\o X)\o Z
      +Q^*(N(\o X),N(\o Y))\o Z\\
      \phantom{\ \qquad\qquad\qquad}=R(\o X,\o Y)\o Z
      +\Omega(\beta\o X,\beta\o Y)\o Z$,\newline
where\newline $\Omega(\beta\o X,\beta\o Y)\o Z
      = (\nabla_{\beta\o Y}B)(\o X,\o Z)
          -(\nabla_{\beta\o X}B)(\o Y,\o Z)+A(R(\o X,\o Y)\o\eta,\o Z)\\
\phantom{..\qquad\qquad\qquad} +B(\o Y,B(\o X,\o Z))-B(\o X,B(\o
Y,\o Z)).$

\item[(b)]$P^*(\o X,\o Y)\o Z+Q^*(\o X,N(\o Y))\o Z=P(\o
X,\o Y)\o Z
      +\Omega(\gamma\o X,\beta\o Y)\o Z$,\newline
where\newline $\Omega(\gamma\o X,\beta\o Y)\o Z
      =-(\nabla_{\gamma\o X}B)(\o Y,\o Z)
          +(\nabla_{\beta\o Y}A)(\o X,\o Z)+A(P(\o X,\o Y)\o\eta,\o Z)\\
\phantom{.\ \qquad\qquad\qquad}-B(T(\o X,\o Y),\o Z)+B(\o Y,A(\o
X,\o Z))
          -A(\o X,B(\o Y,\o Z)).$

  \item[(c)]$Q^*(\o X,\o Y)\o Z=Q(\o X,\o Y)\o Z
      +\Omega(\gamma\o X,\gamma\o Y)\o Z$,\newline
where\newline $\Omega(\gamma\o X,\gamma\o Y)\o Z=(\nabla_{\gamma\o
Y}A)(\o X,\o Z)
      -(\nabla_{\gamma\o X}A)(\o Y,\o Z)+A(\o Y,A(\o X,\o Z))
      -A(\o X,A(\o Y,\o Z))$.
\end{description}
\end{prop}

\begin{cor}\label{cor.2}
\end{cor}
\begin{description}

\item[(a)] Assume that the $\pi$-tensor field $B$ vanishes. Then, $R^*=0$
      if, and only if, $R=0$.

\item[(b)] The $\pi$-tensor field $N$ vanishes if, and only if, $N_0$
vanishes.

\item[(c)] $A(\o X,\o Y)=A(\o Y,\o X)$,
      that is, the $\pi$-tensor field $A$ is symmetric.
\end{description}

\prof  (a) follows from Proposition \ref{prop.3} (a), taking the
fact that $A(\o X,\o\eta)=0$ into account.\newline (b) follows from
Proposition \ref{prop.3}(b) and Proposition \ref{prop.2}(b).\newline
(c) follows from Proposition \ref{prop.3} (c).\quad $\Box$

\begin{lem}\label{lem.3}
For all $\ \o X,\o Y\in\cp$, we have
\begin{description}

\item[(a)] $\
(\nabla_{\gamma\o X}\omega)(\o Y)=-\omega(T(\o X,\o Y))$.

\item[(b)] $\ (\nabla_{\beta\o X}\omega)(\o Y)=\ell^*(B(\o X,\o Y))
+\ell^*(A(N(\o X),\o Y))+L^{-1}h(N(\o X),\o Y)-\omega(T(N(\o X),\o
Y))$.
\par In particular,
$\ (\nabla_{\gamma\o X}\omega)(\o\eta)=0\ $ and $\ (\nabla_{\beta\o
X}\omega)(\o\eta)=\ell^*(N(\o X))$.

\end{description}
\end{lem}

\prof (a) Using Lemma \ref{lem.1} and equations (\ref{eq.5}), we get
$$(\nabla_{\gamma\o X}\ell^*)(\o Y)=L^{-1}h(\o X,\o Y)
+(\nabla_{\gamma\o X}\omega)(\o Y).$$ As $\ \ell(T(\o X,\o Y))=0\ $,
Proposition \ref{prop.2} and equations (\ref{eq.6}), (\ref{eq.7})
and (\ref{eq.5}) give
$$(\nabla_{\gamma\o X}\ell^*)(\o Y)={L^*}^{-1}h^*(\o X,\o Y)
-\omega(T(\o X,\o Y)).$$ Taking (\ref{eq.5}) into account, the
result follows from the above two identities.\newline (b) Using
(\ref{eq.9}), (\ref{eq.5}), Lemma \ref{lem.1} and (a) above, we have
$$(\nabla_{\beta^*\o X}\ell^*)(\o Y)=(\nabla_{\beta\o X}\omega)(\o Y)
-L^{-1}h(N(\o X),\o Y)+\omega(T(N(\o X),\o Y)).$$ On the other hand,
using equations (\ref{eq.6}), (\ref{eq.7}) and (\ref{eq.9}) and
Lemma \ref{lem.1}, we get
$$(\nabla_{\beta^*\o X}\ell^*)(\o Y)=\ell^*(A(N(\o X),\o Y))
+\ell^*(B(\o X,\o Y)).$$ The result follows then from the above two
equations. \quad $\Box$
\par
The next result gives an explicit expression for the $\pi$-tensor
field $A$.

\begin{prop}\label{prop.4}
The $\pi$-tensor field $A$ is given by\:
$$ A = \frac1{2L^*}(h\otimes\o m +\nu\otimes\phi+\phi\otimes\nu)
             -\frac1{2{L^*}^2}\{2L^*(\omega\otimes\o\eta)T
             +2\nu\otimes\nu\otimes\o\eta+\nu(\o m)h\otimes\o\eta\}.$$
\end{prop}

\prof  We first prove, for all $\o X,\o Y,\o Z\in\cp$, that

\begin{eqnarray}\label{eq.10}
A^*(\o X,\o Y,\o Z) &=& \frac1{2L}\{h(\o X,\o Y)\nu(\o Z)
                        +h(\o Y,\o Z)\nu(\o X)
                        +h(\o X,\o Z)\nu(\o Y)\}\nonumber\\
                    & & \mbox{}-\omega(T(\o X,\o Y))\ell^*(\o Z).
\end{eqnarray}
Using equations (\ref{eq.6}) and (\ref{eq.7}), taking the fact that
$\nabla^*g^*=0$ into account, we have
$$(\nabla_{\gamma\o X}g^*)(\o Y,\o Z)=A^*(\o X,\o Y,\o Z)
+A^*(\o X,\o Z,\o Y).$$ On the other hand, using (\ref{eq.5}), Lemma
\ref{lem.1} and Proposition \ref{prop.1}, we get
\begin{eqnarray*}
(\nabla_{\gamma\o X}g^*)(\o Y,\o Z) &=& L^{-1}\{h(\o X,\o Y)\nu(\o
Z)
+h(\o Y,\o Z)\nu(\o X)+h(\o X,\o Z)\nu(\o Y)\}\\
& &+\ell^*(\o Y)(\nabla_{\gamma\o X}\omega)(\o Z) +\ell^*(\o
Z)(\nabla_{\gamma\o X}\omega)(\o Y).
\end{eqnarray*}
Taking Corollary \ref{cor.1} and Lemma \ref{lem.3}(a) into account,
(\ref{eq.10}) follows from the above two equations.
\par
Now, using equations (\ref{eq.5}), (\ref{eq.3}) and (\ref{eq.2}),
taking Proposition \ref{prop.1} into account, one can show that:
\begin{equation}\label{eq.11}
  \left.
      \begin{array}{rcl}
         g(\o m,\o Z) &=& \tau^{-1}g^*(\phi^*(\o m),\o Z)\\
         h(\o X,\o Z) &=& \tau^{-1}g^*(\phi^*(\o X),\o Z)
      \end{array}
  \right\}
\end{equation}
\par
Substituting (\ref{eq.11}) into (\ref{eq.10}), taking (\ref{eq.2})
into account, it follows from the nondegeneracy of $g^*$ that
$$A(\o X,\o Y)=\frac1{2L^*}\{h(\o X,\o Y)\phi^*(\o m)
+\nu(\o X)\phi^*(\o Y)+\nu(\o Y)\phi^*(\o X)\}
-\frac1{L^*}\omega(T(\o X,\o Y))\o \eta.$$
\par The result follows then from the fact that
$\phi^*=\phi-{L^*}^{-1}\nu\otimes\o\eta.$ \quad $\Box$

\begin{cor}\label{cor.3}
For all $\o X,\o Y,\o Z\in\cp$, we have\:
$$T^*(\o X,\o Y,\o Z)=\tau T(\o X,\o Y,\o Z)
+\frac1{2L}\{h(\o X,\o Y)\nu(\o Z) +h(\o Y,\o Z)\nu(\o X)+h(\o X,\o
Z)\nu(\o Y)\}.$$
\par {\em In fact, this formula follows from Proposition \ref{prop.2}(c)
and from equation (\ref{eq.10}).}
\end{cor}

\smallskip
\par
Let $(x^i),\ i=1,2,\ldots,n$, be a system of local coordinates on
$M$ and let $(x^i,y^i)$ be the associated canonical system of local
coordinates on $TM$ and $\tm$. The natural bases of $T_u(\tm)$ and
$H_u(\tm)$ are denoted respectively by
$(\partial_i,\partial_{,i})_u$ and $(e_i)_u$. The values of the lift
$\o\partial_i$ of $\partial_i$ at $u$ form a basis for the fibre
over $u$ in $\p$. We write
$$\nabla_{\partial_i}\,\o\partial_j=\Gamma_{ij}^h\,\o\partial_h,\quad
\nabla_{\partial_{,i}}\,\o\partial_j=C_{ij}^h\,\o\partial_h,\quad
\nabla_{e_i}\,\o\partial_j=\o\Gamma_{ij}^h\,\o\partial_h.$$ The
relations (\ref{eq.5}) can be expressed locally by:

\begin{equation}\label{eq.12}
  \left.
     \begin{array}{rcl}
{\ell^*}_i   &=& \ell_i+b_i\\
{h^*}_{ij}   &=& \tau h_{ij}\\
{g^*}_{ij}   &=& \tau(g_{ij}-\ell_i\ell_j)+{\ell^*}_i{\ell^*}_j\\
{g^*}^{ij}   &=&
\tau^{-1}g^{ij}+\mu\ell^i\ell^j-\tau^{-2}(\ell^ib^j+\ell^jb^i)
     \end{array}
  \right\}
\end{equation}
where $\ \ell^i=y^i/L,\ \ell_i=g_{ir}\ell^r,\
\mu=(Lb^2+\alpha)/L^*\tau^2, \ b^2=g_{ij}b^i(x)b^j(x)=b_ib^i\ $;
$b^i(x)$ being the components of the $\pi$-vector field $\o b$. We
use the notations:

\begin{equation}\label{eq.13}
\left.
\begin{array}{rclcl}
b_{ij}   &:=& \nabla_{e_i}b_j,\qquad\qquad\quad b_{i0}   &:=& b_{ik}y^k\\
b_{[ij]} &:=& (b_{ij}-b_{ji})/2,\qquad\         b_{(ij)} &:=&
(b_{ij}+b_{ji})/2
\end{array}
\right\}
\end{equation}
\par
After some lengthy but straightforward calculations, using
(\ref{eq.12}) and (\ref{eq.13}), the $\pi$-tensor fields $N_0$, $N$
and $B$ in (\ref{eq.7}) are given locally by:
\begin{eqnarray*}
N_0^h    &=& {\Gamma^*}_{00}^h-\Gamma_{00}^h\\
         &=& {\ell^*}^hb_{(00)}-2L^*{g^*}^{hr}b_{[r0]}.\\
N_i^h    &=& {\Gamma^*}_{i0}^h-\Gamma_{i0}^h\\
         &=& {g^*}^{hk}\left(L^*b_{[ik]}-{\ell^*}_i b_{[k0]}\right)
             +{\ell^*}^hb_{(i0)}+(1/{2L^*}){h^*}_i^hb_{00}
             +2L^*{g^*}^{rk}{C^*}_{ir}^hb_{[k0]}.
\end{eqnarray*}
\begin{eqnarray}\label{eq.14}
B_{ij}^h &=& {\Gamma^*}_{ij}^h-\o\Gamma_{ij}^h-\Gamma_{i0}^r{C^*}_{jr}^h\nonumber\\
         &=& {g^*}^{hr}\left({\ell^*}_ib_{[jr]}+{\ell^*}_jb_{[ir]}\right)
             +{\ell^*}^hb_{(ij)}\nonumber\\
         & & {}+(1/{2L^*})\left(b_{i0}{h^*}_j^h
             +b_{j0}{h^*}_i^h-{g^*}^{hk}b_{k0}{h^*}_{ij}\right)\nonumber\\
         & & {}+{g^*}^{hp}{C^*}_{ijr}
             \left\{{g^*}^{rk}\left(L^*b_{[pk]}
             -{\ell^*}_pb_{[k0]}\right)+{\ell^*}^rb_{(p0)}
             +(1/{2L^*})b_{00}{h^*}_p^r\right.\nonumber\\
         & & \left.{}+2L^*{g^*}^{mk}{C^*}_{pm}^rb_{[k0]}\right\}
             -{g^*}^{hp}{C^*}_{irp}
             \left\{{g^*}^{rk}\left(L^*b_{[jk]}
             -{\ell^*}_jb_{[k0]}\right)\right.\nonumber\\
         & & \left.{}+{\ell^*}^rb_{(j0)}
             +(1/{2L^*})b_{00}{h^*}_j^r
             +2L^*{g^*}^{mk}{C^*}_{jm}^rb_{[k0]}\right\},
\end{eqnarray}
where $\ h_i^r=g^{rj}h_{ij}$.
\par
Moreover, Lemma \ref{lem.3}(b) may be expressed locally by:

\begin{equation}\label{eq.15}
b_{ij}=B_{ij}^k{\ell^*}_k+A_{rj}^kN_i^r{\ell^*}_k
+L^{-1}N_i^rh_{rj}-b_kN_i^rT_{rj}^k.
\end{equation}

\Section{Main Results}
\par
Let $(M,L^*)$ be a generalized Randers manifold with $(M,L)$ as its
associated Finsler manifold.
\par
Using formulae (\ref{eq.14}) and (\ref{eq.15}), one can prove

\begin{thm}\label{thm.1}
The $\pi$-tensor field $B$ vanishes if, and only if, the
$\nabla$-horizontal \linebreak covariant derivative of $\omega$
vanishes {\em(i.e. $\nabla_{\beta \o X}\omega=0 \quad \forall \o
X$)}.
\end{thm}

\par
Let $J$ be the vector $1$-form on $\tm$ defined by
$J=\gamma\circ\rho$, then $d_J\alpha=d\alpha\circ J$, where $\alpha$
is the function defined by (\ref{eq.4}). The proof of the following
result is similar to that of Proposition 2 of ~\cite{[11]}.

\begin{lem}\label{lem.4}
The generalized Randers manifold $(M,L^*)$ and its associated
Finsler manifold $(M,L)$ have both the same geodesics if, and only
if, $d_J\alpha$ is closed.
\end{lem}

\begin{thm}\label{thm.2}
For the generalized Randers manifold $(M,L^*)$ and its associated
Finsler manifold $(M,L)$, the following assertions are equivalent
\begin{description}
    \item[(a)] $(M,L)$ and $(M,L^*)$ have the same
geodesics.

\item[(b)] $B(\o V,\o V)$ vanishes for all $\o V$, where $\o
V=\o\eta\mid_{\ti c(t)}$.

\item[(c)] $d_J\alpha$ is closed.

 \item[(d)] $N$ vanishes identically.

\end{description}
\end{thm}

\prof (a)$\Longleftrightarrow$ (b): Theorem 2 of
~\cite{[12]}.\newline (b)$\Longleftrightarrow$ (c): Lemma
\ref{lem.4} and Theorem 2 of ~\cite{[12]}.\newline
(c)$\Longleftrightarrow$ (d): Lemma \ref{lem.4}, Corollary
\ref{cor.2}(b) and Theorem 2 of ~\cite{[12]}. \quad $\Box$
\medskip
\par
Corollary \ref{cor.2}(a) and Theorem \ref{thm.1} imply

\begin{prop}\label{prop.5}
Let $\ \nabla_{\beta\o X}\omega=0\ $ for all $\o X\in\cp$. Then,
$R^*$ vanishes if, and only if, $R$ vanishes.
\end{prop}

\par
A Finsler manifold $(M,L)$ is a Berwald manifold ~\cite{[4]} if the
torsion tensor $T$ satisfies the condition that $ \nabla_{\beta\o
X}T=0\ $ for every $\o X\in\cp$. A Finsler manifold $(M,L)$ is
locally Minkowskian ~\cite{[4]} if, and only if, $R=0$ and
$\nabla_{\beta\o X}T=0$.
\par
Combining Theorems 5 and 6 of ~\cite{[12]} and Theorem \ref{thm.1},
we get
\begin{thm}\label{thm.3}
For the generalized Randers manifold $(M,L^*)$ and its associated
Finsler manifold $(M,L)$, suppose that $\ \nabla_{\beta\o X}\omega=0
\quad\forall\o X\in\cp$.\newline Let $(M,L)$ {\em(resp. $(M,L^*)$)}
be a Berwald {\em (or locally Minkowskian)} manifold. A necessary
and sufficient condition for $(M,L^*)$ {\em(resp. $(M,L)$)} to be a
Berwald {\em(or locally Minkowskian)} manifold is that $\
\nabla_{\beta\o X}A=0\ $ for all $\ \o X\in\cp$.
\end{thm}

\par
A Finsler manifold $(M,L)$ is a Landsberg manifold ~\cite{[5]} if it
satisfies the condition that $\ P(\o X,\o Y)\o\eta=0\ $ for all $\o
X,\o Y\in\cp$.
\par
Combining Theorem 7 of ~\cite{[12]} and Theorem \ref{thm.1}, we get

\begin{thm}\label{thm.4}
For the generalized Randers manifold $(M,L^*)$ and its associated
Finsler manifold $(M,L)$, suppose that $\nabla_{\beta\o X}\omega=0
\quad\forall\o X\in\cp$.\newline $(M,L^*)$ is a Landsberg manifold
if, and only if, $(M,L)$ is a Landsberg manifold.
\end{thm}
\par
Combining Theorem 1 of ~\cite{[12]} and Theorem \ref{thm.2}, we get

\begin{thm}\label{thm.5}
For the generalized Randers manifold $(M,L^*)$ and its associated
Finsler manifold $(M,L)$, suppose that the $1$-form $d_J\alpha$ is
closed. The horizontal distribution of $(M,L^*)$ is completely
integrable if, and only if, the horizontal distribution of $(M,L)$
is completely integrable.
\end{thm}
\par
A general Landsberg manifold ~\cite{[5]} is a Finsler manifold such
that the trace of the linear map $\ \o Y\longmapsto P(\o X,\o
Y)\o\eta\ $ is zero, for all $\pi$-vector fields $\o X$. It is
characterized by the condition that $\ \nabla_{\beta\o\eta}C=0$,
where $C$ is the $\pi$-form obtained from the torsion tensor $T$ by
contraction.

\begin{lem}\label{lem.5}
The $\pi$-forms $C$ and $C^*$ are related by\:
$$C^*=C+\frac{n+1}{2L^*}\,\nu.$$
\end{lem}
\par
The proof of this lemma is similar to that found in ~\cite{[10]}.

\begin{prop}\label{prop.6}
For all $\pi$-vector field $\o X\in\cp$, we have\:
\begin{eqnarray*}
({\nabla^*}_{\beta^*\o\eta}C^*)(\o X) &=&
(\nabla_{\beta\o\eta}C)(\o X)
     -(\nabla_{\gamma N_0}C)(\o X)+C(A(N_0,\o X))\\
& &  {}-C(B(\o\eta,\o X))+\frac{n+1}{2L^*}
     \left\{(\nabla_{\beta\o\eta}\nu)(\o X)
     -(\nabla_{\gamma N_0}\nu)(\o X)\right.\\
& &  \left. {}+\nu(A(N_0,\o X))-\nu(B(\o\eta,\o X))\right\}.
\end{eqnarray*}
\par
In particular, if the $1$-form $d_J\alpha$ is closed, then
$${\nabla^*}_{\beta^*\o\eta}C^*=\nabla_{\beta\o\eta}C
+\frac{n+1}{2L^*}\,\nabla_{\beta\o\eta}\nu.$$
\end{prop}

\prof The first formula follows from Lemma \ref{lem.5}. The second
formula follows from the first one, Theorem \ref{thm.2} and from
Corollary \ref{cor.2}(b).\quad $\Box$

\medskip
\par
Now, Proposition \ref{prop.6} implies
\begin{thm}\label{thm.6}
For the generalized Randers manifold $(M,L^*)$ and its associated
Finsler manifold $(M,L)$, suppose that $d_J\alpha$ is
closed.\newline Let $(M,L)$ {\em(resp. $(M,L^*)$)} be a general
Landsberg manifold. A necessary and sufficient condition for
$(M,L^*)$ {\em(resp. $(M,L)$)} to be a a general Landsberg manifold
is that $\nabla_{\beta\o\eta}\nu=0$.
\end{thm}
\par
Finally, by Theorems 3 and 4 of ~\cite{[12]} and Theorem
\ref{thm.2}, we have

\begin{thm}\label{thm.7}
If $d_J\alpha$ is closed, the geodesics of the generalized Randers
manifold and those of its associated Finsler manifold have both the
same Morse index.
\end{thm}

\end{document}